\documentclass[final,3p]{elsarticle}
\usepackage[english]{babel}
\usepackage[utf8]{inputenc}
\usepackage{typearea}
\usepackage{latexsym}
\usepackage{lineno,hyperref}
\usepackage{amsmath,amssymb,amsthm,mathtools}
\usepackage{algorithm}
\usepackage{algorithmic}
\usepackage{graphicx}
\usepackage[abs]{overpic}

\newcommand{\R}{\mathbb{R}}

\newcommand{\CC}{\mathcal{C}}
\def\M{\mathcal{M}}

\newcommand{\Sx}{\mathcal{S}}
\newcommand{\Dx}{\mathcal{D}}

\newcommand{\df}{Diff^+}

\newcommand{\vphi}{\varphi}

\newcommand*{\quark}{\setbox0\hbox{$x$}\hbox to\wd0{\hss$\cdot$\hss}}

\usepackage{xcolor}
\newcommand{\todo}[1]{\bgroup\color{red}#1\egroup}

\newtheorem{proposition}{Proposition}
\newtheorem{lemma}{Lemma}
\newtheorem{example}{Example}
\begin{document}
\begin{frontmatter}
\title{
On Geodesics in the Spaces of Constrained Curves
}
\author[first]{Esfandiar Nava-Yazdani\corref{cor}} \cortext[cor]{Corresponding author} \ead{navayazdani@zib.de}
\begin{abstract}In this work, we study the geodesics of the  space of certain geometrically and physically motivated subspaces of the space of immersed curves endowed with a first order Sobolev metric. This includes elastic curves and also an extension of some results on planar concentric circles to surfaces. The work focuses on intrinsic and constructive approaches.
\end{abstract}
\begin{keyword}
 Sobolev metric \sep Riemannian submersion \sep Shape analysis \sep Elastic curve \sep Longitudinal analysis 
\end{keyword}
\end{frontmatter}

\section{Introduction}
The study of curves and their shapes is an emerging research field with numerous and varied application areas such as computer vision, image analysis and morphology, and has attracted a great deal of attention over the past years. While the study of planar closed curves advances the approaches to 2D shapes, many applications naturally lead to manifold-valued curves. Some prominent examples are provided by curves in Lie groups such as the Euclidean motion group, or more generally, symmetric spaces including the Grassmannian and the Hadamard-Cartan manifold of positive definite matrices. Furthermore, an essential task in analysis of longitudinal data, is to compare curves in a high or infinite dimensional space. These applications motivate the study of reparametrization invariant metrics on spaces of curves and their shapes. A Riemannian framework for analysis in these spaces is desirable, because it naturally provides these spaces with a rich structure for whose treatment powerful tools are available. 

Remarkably, it has been shown in~\cite{michor2005vanishing} that the simple natural candidate, the $L^2$-metric, always vanishes. This has motivated the investigation of stronger Sobolev metrics. The space of Euclidean curves under these metrics has been widely studied. Particular, the works ~\cite{michor2006riemannian,michor2007overview,bruveris2014geodesic,bauer2016use} address some core properties including (geodesic and metric) completeness and the geodesic equation. Some of these results have been extended to manifold-valued curves in~\cite{bauer2020SobolevMO}. In general, the numerical realization of the underlying infinite dimensional Riemannian calculus, concerning both the space of curves and the reparametrization group acting on it, poses enormous computational challenges. Several numerical approaches to computation of geodesics and distances for Euclidean data have been proposed over the past years. In particular, we refer to the overview~\cite{bauer2020intrinsic} and~\cite{bauer2017varifold, bauer17numcurve, bauer18relaxed}. Separate from the active research at the intersection of structured data and shape analysis, there is an evolving body of work in geometric statistics demonstrating how exploiting the intrinsic geometric structure of sequential data yields effective practical tools for relevant statistics and analysis tasks. To name a few works,~\cite{sri2010protein, su2014statistical,sri2016functional, bauer2015med, nava2020geo, hanik2020} and ~\cite{nava2023hur}, employed a Riemannian approach to provide certain summary statistics (mean and principle modes) for real world applications like analysis of protein structures, bird migration, HeLa cell nuclei, osteoarthritis, deformations during cardiac cycles and hurricane tracks. 

This paper is organized as follows. Section \ref{sec:1} presents the Riemannian setting and notations. Section \ref{sec:2} is devoted to applications on paths of unparametrized curves. Therein, we present the variation formulae and a result on conservation of the curvature, an extension of a result from~\cite{michor2006riemannian} on plane circles to arbitrary surfaces, and as further example of geodesics in the space of constrained curves, shortest paths of elastica. 
\section{Riemannian Framework}\label{sec:1}
Let $(M,g )$ be a finite dimensional Riemannian manifold and $\M$ the Fr{\'e}chet manifold of smooth immersed curves from $\Dx$ in $M$, where $\Dx$ denotes either the unit circle $S^1$ or the unit interval $I$ for closed or open curves, respectively. Tangent space of $\M$ at a curve $c$ is the space of vector fields along $c$. Let us denote the group of orientation preserving diffeomorphisms on $\Dx$ by $\df$. We are interested in a reparametrization invariant Riemannian metric $G$ on $\M$: 
\[
G_{c\circ \vphi}(h\circ \vphi,k\circ\vphi)=G_c(h,k),
\]
for any $c\in \M$,\, $h,k\in T_c\M$ and $\vphi\in\df$. The above equivariance guarantees that the induced distance satisfies 
\[
d(c_0\circ\vphi,c_1\circ\vphi)=d(c_0,c_1),
\]
for any two curves $c_0$ and $c_1$. The induced distance on the space of unparametrized curves
\[
\Sx=\M/\df
\]
reads
\begin{align*}
d^S([c_0],[c_1])&=inf\,\{d(c_0, c_1\circ\vphi):\, \vphi\in \df\}\\
&=inf\,\{d( c_0\circ\vphi,c_1):\, \vphi\in \df\}.
\end{align*}
Denoting the quotient map by $\pi$, we have the canonical decomposition of the tangent bundle $T\M$ into vertical subbundle $Ver=ker(d\pi)$ and its orthogonal complement, the horizontal subbundle $Hor$. More explicitly, denoting $\CC=\{a\in C^\infty(\Dx,\R):\, a(0)=a(1)=0 \text{ if }\Dx=I\}$,
\begin{align*}
Ver_c&=\{a\dot{c}:\, a\in \CC\},\\
Hor_c&=\{h\in T_c\M:\, G_c(h,a\dot{c})=0\text{ f. a. }a\in \CC\},
\end{align*}
where dot stands for differentiation with respect to $t$.
We recall, that the quotient map $\pi$ becomes a Riemannian submersion (a comprehensive discussion is given in~\cite{michor2007overview}), thus encodes the geometry of $\Sx$ by horizontal lifting. This fact and the above decomposition is also consensus in other shape theories, including  the landmark-based approach of Kendall. Now, let $\nabla$ be the Levi-Civita connection of $M$ and denote $\omega = |\dot{c}(t)|$ $d\theta=\omega\, dt$, and the unit tangent field of $c$ by $T$.  Remarkably, it has been shown in~\cite{michor2005vanishing} that the simple natural candidate, the $L^2$-metric
\[
G_c^{L^2}(h,k)=\int_\Dx g_c( h,k) \,d\theta
\]
always vanishes and cannot be used to distinguish the shapes of curves. This has motivated the investigation of higher order Sobolev metrics. In this work, we consider the first order one, given by 
\[
G_c(h,k)= \int_\Dx g_c( h, k)+g_c( \nabla_Th, \nabla_Tk)\, d\theta,
\]
which is well suited for numerical computations and among the most widely used metrics in applications.

Let $C$ be a smooth path with values in $\M$. For brevity, we frequently write $c(s,t):=(C(s))(t)$ with $s\in I:=[0,1]$ and $t\in\Dx$, and denote the speed of $C$ by $\nu$, i.e., $\nu =\sqrt{G_c(c^\prime,c^\prime)}$. Here, prime stands for differentiation with respect to the path parameter $s$. Thus, the geodesic differential equation for $C$ in integrated form reads $\nu^\prime =0$. For $t\mapsto c(.,t)$, we denote the signed curvature by $\kappa$ and differentiation with respect to arc length using subscript as in $\kappa_\theta=\partial_\theta\kappa$. Thus, $\nabla_TT=\kappa N$ is the curvature vector field and $N$ a unit normal field along $t\mapsto c(.,t)$ on $I$. We recall that $C$ is called horizontal, iff its tangent field $C^\prime$ is horizontal. Denoting 
$\eta= c^\prime-\nabla^2_Tc^\prime$, integration by parts shows that $C$ is horizontal iff $g(\eta,T)=0$. We denote the normal component of $C^\prime$ by $\rho$, i.e., $\rho=g(c^\prime, N)$. Important special cases for $c$, often canonical in geometric variational tasks, are normality, i.e., $g(c^\prime,T)=0$ (clearly, equivalent to $L^2$-horizontality) and particularly $c^\prime=\rho N$.
\section{Applications}\label{sec:2}
\subsection{Variation formulae and Conserved Quantities}
Next, we consider the variations of the main scalar geometric quantities, the length element $\omega$ and the curvature $\kappa$ of the $t$-curves along $C$.
\begin{lemma}\label{lem1}
Let $K$ denote the sectional curvature. Then, the following holds.\\
a) 
\begin{align*}
\omega^\prime &= g(\nabla_Tc^\prime,T)\omega,\\
\kappa^\prime &= g(\nabla^2_Tc^\prime,N)-2\kappa g(\nabla_Tc^\prime,T)+K(N,T)g(c^\prime,N).
\end{align*}
b) 
Suppose that $C$ is normal, i.e., $c^\prime=\rho N$. Then, we have
\begin{align*}
\omega^\prime = -\rho\kappa\omega.
\end{align*}
Moreover, $C$ is horizontal iff $\rho^2\kappa$ is constant along the $t$-curves, i.e., $(\rho^2\kappa)_\theta =0$.
\begin{proof}
a) Derivation of the first identity is straightforward (cf.~\cite{bauer2020SobolevMO}). Let $R$ denote the curvature tensor and $[.\,,\,.]$ the Lie bracket. Due to $\kappa=g(\nabla_TT,N)$ we have
\begin{align*}
    \kappa^\prime &= g(\nabla_{c^\prime}\nabla_TT, N)+g(\nabla_TT,\nabla_{c^\prime}N)\\
    & = g(\nabla_{c^\prime}\nabla_TT, N)+\kappa g(N,\nabla_{c^\prime}N)\\
    &=g(\nabla_{c^\prime}\nabla_TT, N)\text{ (since $g(N,N)=1$)}\\
    &=g(\nabla_T\nabla_{c^\prime}T,N)+g(\nabla_{[{c^\prime},T]}T,N)+g(R({c^\prime},T)T,N).
\end{align*}
Due to $\nabla_{c^\prime}T=\nabla_Tc^\prime-\frac{\omega^\prime}{\omega}T$, we have $[c^\prime,T]=-g(\nabla_Tc^\prime,T)T$ and the first term in the sum can be written as
\begin{align*}
    g(\nabla_T\nabla_{c^\prime}T,N)&=g(\nabla_T(\nabla_Tc^\prime-g(\nabla_Tc^\prime,T)T),N)\\
    &=g(\nabla^2_T{c^\prime},N)- g(\nabla_T{c^\prime},T)g(\nabla_TT,N)\\
    &=g(\nabla^2_T{c^\prime},N)-\kappa g(\nabla_T{c^\prime},T).
\end{align*}
Furthermore, the second term can be written as
\begin{align*}g(\nabla_{[c^\prime,T]}T,N)&=- g(\nabla_Tc^\prime,T)g(\nabla_TT,N)\\ &=-\kappa g(\nabla_Tc^\prime,T).
\end{align*}
b) We have 
\begin{align*}
\omega^\prime &= g(\nabla_T(\rho N),T)\omega= g(\rho_\theta N+\rho\nabla_TN,T)\omega\\ &=-\rho g(N,\nabla_TT)\omega=-\rho\kappa\omega.
\end{align*}
Furthermore, $g(c^\prime, T)=0$ implies
\begin{align*}
    g(\eta, T)&= g(\nabla_T^2c^\prime, T)\\ &= g(\rho_{\theta\theta}N+2\rho_\theta \nabla_TN+\rho \nabla_T^2N, T)\\
              &= -2\rho_\theta\kappa -\rho g(\nabla_T(\kappa T),T)\\ &=-2\rho_\theta\kappa-\rho\kappa_\theta.
\end{align*}
As $\rho$ does not vanish, multiplying the above identity with $\rho$, we conclude that, $C$ is horizontal, i.e., $g(\eta, T)$ disappears iff $(\rho^2\kappa)_\theta=0$.
\end{proof}
\end{lemma}
\begin{proposition}
    Suppose that $dim(M)=2$ and $C$ is normal and horizontal. Then, the following holds.\\
    a) The unit tangent field $T$ is parallel along $C$ iff $\kappa_\theta=0$, i.e., all $t$-curves have constant curvature. In this case, $C$ is a geodesic iff
    $(\omega^\prime)^2=\frac{\alpha\omega\kappa^2}{1+\kappa^2}$ with a positive constant $\alpha$.\\
    b) $\kappa$ is constant along $C$, i.e., $\kappa^\prime=0$ iff $2\kappa\kappa_{\theta\theta}=3\kappa_{\theta}^2+4\kappa^2(\kappa^2+K)$.
    \begin{proof}
    a) We have
    \begin{align*}
    \nabla_{c^\prime}T &=g(\nabla_{c^\prime}T,N)N\\
    & = g(\nabla_Tc^\prime-\frac{\omega^\prime}{\omega}T,N)N\\
    &=\rho_\theta N
    \end{align*}
    Hence, $T$ is parallel along $C$ iff $\rho_\theta=0$, which is due to
    horizontality and lemma \ref{lem1} equivalent to $\kappa_\theta=0$. Thus, denoting the lengths of the $t$-curves by $L$, we have $\nu^2=\rho^2(1+\kappa^2)L$. Hence, $\nu^\prime=0$ iff $(\rho^2(1+\kappa^2)L)^\prime=0$. Utilizing $\omega^\prime=-\rho\kappa\omega$, we arrive at the desired equation.\\
    b) Lemma \ref{lem1} implies
    \[
    \kappa^\prime = \rho_{\theta\theta}+\rho(\kappa^2 + K).
    \]
     Moreover, due to horizontality, we have $2\rho_\theta\kappa+\rho\kappa_\theta=0$ implying $2\rho_{\theta\theta}\kappa+3\rho_\theta\kappa_\theta+\rho\kappa_{\theta\theta}=0$. Therefore, $4\rho_{\theta\theta}\kappa^2-3\rho\kappa_\theta^2+2\rho\kappa\kappa_{\theta\theta}=0$. Inserting in the expression for $\kappa^\prime$, completes the proof. 
    \end{proof}
\end{proposition}
The following is an immediate applications of lemma \ref{lem1}.
\begin{example}
    An integral curve of the shortening flow given by $c^\prime=\kappa N$ is horizontal iff $\kappa_\theta=0$.
\end{example}
\subsection{Concentric Circles}
In the following, we extend a result from~\cite{michor2006riemannian}) on the horizontality of a family of concentric plane circles and the corresponding geodesic equation to arbitrary surfaces.
\begin{proposition}
    Let $dim(M)=2$ and $I\ni s\mapsto c(s,.)\in\M$ a family of concentric circles with radius $r$. Then, $c$ is horizontal iff the circles have constant curvature. Moreover, the following holds.\\
    a) $C$ is a geodesic iff 
    \begin{equation}
        \left((r^\prime)^2\int_0^{2\pi}\omega + \frac{(\omega_r)^2}{\omega}\,dt\right)^\prime =0.\label{eq:GDE}
    \end{equation}
    b) If $M$ has constant curvature, then $C$ is horizontal and the geodesic differential equation reads
    \begin{equation}
        \left((r^\prime)^2 \left(\omega + \frac{(\omega_r)^2}{\omega} \right)\right)^\prime =0,\label{eq:GDE_circ}
    \end{equation}
    with
    \[
    \omega (r)= \begin{cases}
    \frac{1}{\sqrt{K}}\sin(\sqrt{K} r)\quad & K >0,\\
    r\quad & K =0,\\
    \frac{1}{\sqrt{-K}}\sinh(\sqrt{-K} r)\quad & K < 0.
    \end{cases}
    \]
\begin{proof}
a) Let $p$ denote the center of the circles. Then, we may write $c(s,t)=\exp_p(r\cos t\, e_1+\sin t\, e_2)$, $g$ is given by $dr^2+\omega^2(r,t)dt^2$ and $c^\prime=\rho N$. Liouville's theorem implies $\omega^\prime=r^\prime \kappa \omega$. Due to lemma \ref{lem1}, we have $\omega^\prime=-\rho \kappa \omega$. Thus, $\rho =-r^\prime$. Particularly $\rho_\theta=0$ and due to lemma \ref{lem1}, $c$ is horizontal iff $\kappa_\theta=0$. Thus, the geodesic differential equation, constancy of the speed $\nu$, is given by \eqref{eq:GDE}.\\
b) Now, suppose that $M$ has constant curvature $K$. Then, the circles have constant curvature, i.e. $\kappa_\theta =0$ (cf. proposition 6.11 of~\cite{Lee1997RiemannianMA}). Thus, $c$ is horizontal. The explicit expression for $\omega$ follows immediately from the Jacobi's equation $\omega_{rr}+K\omega=0$ (cf. \cite{docarmo} for details). In particular, $\omega_\theta =0$ and the geodesic equation \eqref{eq:GDE} reduces to \eqref{eq:GDE_circ}.
\end{proof}
\end{proposition}
We remark that, if $M$ is the unit 2-sphere, then the geodesic differential equation \eqref{eq:GDE_circ} is the well-known pendulum equation $\frac{(u^\prime)^2}{\cos u}=const.$, where $u=\frac{\pi}{2}-r$.
\subsection{Elastica}
\begin{figure}[h]
\centering
\includegraphics[width=.23\textwidth]{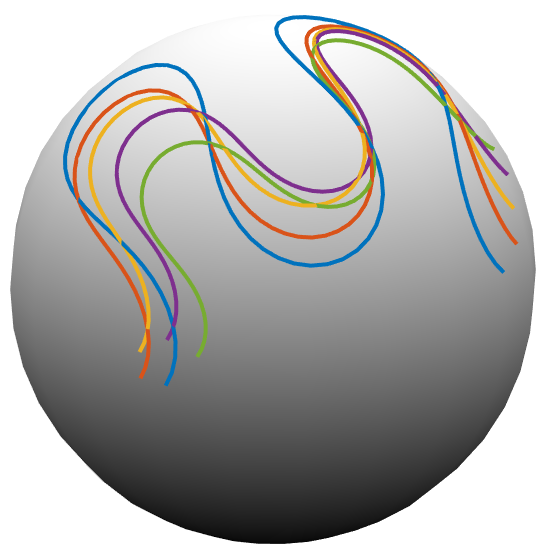}
\hspace{1ex}
\includegraphics[width=.23\textwidth]{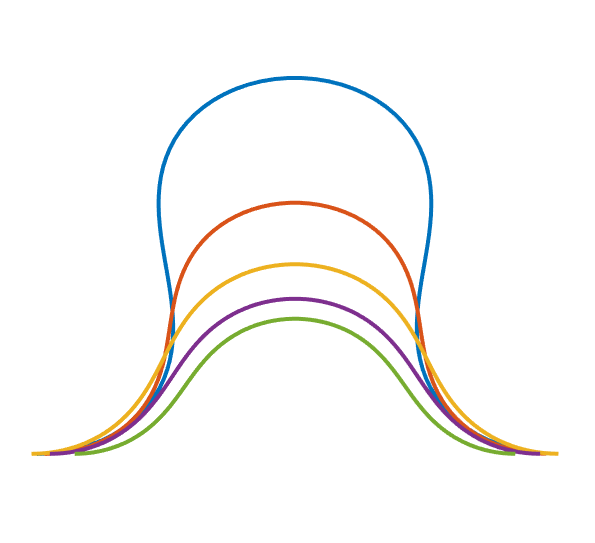}
\hspace{1ex}
\includegraphics[width=.21\textwidth]{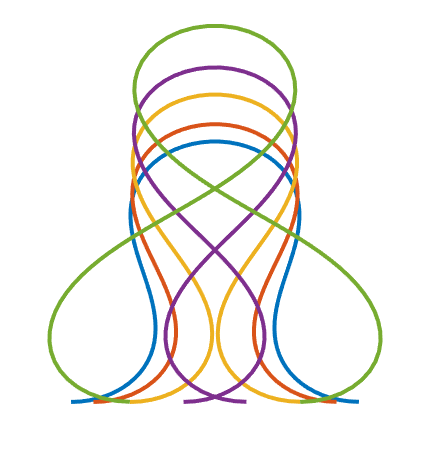}
\hspace{1ex}
\includegraphics[width=.21\textwidth]{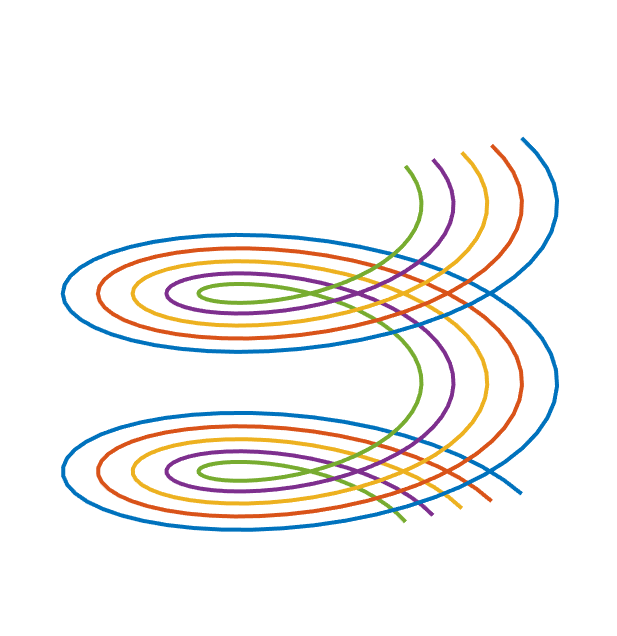}
\caption{Geodesic paths of elastica, joining the green and blue curves, gained by minimizing the energy over the amplitude of curvature $k=\kappa_{max}$, tension $\sigma$ and $\mu=\kappa^2\tau$. These parameters determine the shape of an elastic curve.}
\label{fig:elastica}
\end{figure}

Elastica are minimizers of $\int_D\kappa^2+\lambda\,d\theta$ satisfying given first order boundary data. The Lagrange mulitplier $\lambda$ is the tension. In the following, we assume that the sectional curvature of $M$ is constant. In this case, an elastic curve is also characterized by $\kappa^2\tau = \mu$ for a constant $\mu$ with $\kappa$ given by a Jacobi elliptic function, i.e., $\kappa_{\theta\theta}$ is a cubic polynomial in $\kappa$ (cf.~\cite{Singer2008LecturesOE}). Tension $\lambda$, amplitude (maximum value) of curvature denoted by $k$ and $\mu$ determine curvature and torsion, thus also the shape of the curve (the mentioned examples of circles are paths of elastic curves with $k^6+(2K-\lambda)k^4+2\mu^2=0$, $\kappa=k$). In particular, energy optimization for a shortest paths of elastica, reduces to minimization over these parameters. Figure \ref{fig:elastica} shows some examples.
For $M=\R^3$, coaxial helices with constant pitch constitute a horizontal path of elastica. Moreover, a straightforward computation shows that the following holds.  
\begin{example}
The path of helices given by
\[
c(.,t)=(r\cos t, r\sin t, ht)
\]
is a horizontal geodesic iff $h^\prime = 0$ and 
\[
((r^\prime)^2(\sqrt{r^2+h^2}+\frac{1}{\sqrt{r^2+h^2}}))^\prime = 0.
\]
\end{example}
 We remark that elastic curves can be used to construct Willmore tori (cf.~\cite{Pinkall2005HopfTI}). Thus, as an application, one can construct paths of Willmore tori and investigate their total squared mean curvature.
\section{Conclusion}
In this work, we studied horizontal and geodesic paths in the spaces of manifold-valued curves endowed with reparametrization invariant first order Sobolev metric. Particularly, we presented variation formulas for the curvature of the curves and focusing on geometrically or physically motivated constraints, special cases, for which the geodesic PDE reduces to an ODE. We also presented examples of geodesic paths in the space of elastica, for which geodesic optimization simplifies.  
\section{Acknowledgment}
This work was supported through the German Research Foundation (DFG) via individual funding (project ID 499571814).
\bibliography{main}
\end{document}